\documentclass[a4paper,11pt]{amsart}

\usepackage{amsmath,amsthm,amssymb}
\usepackage[mathscr]{eucal}
\usepackage{cite}
\usepackage[bookmarks,bookmarksnumbered]{hyperref}

\numberwithin{equation}{section}

\theoremstyle{plain}
\newtheorem{theorem}{Theorem}[section]

\newtheorem{lemma}[theorem]{Lemma}

\theoremstyle{definition}
\newtheorem{definition}[theorem]{Definition}

\newtheorem{remark}[theorem]{Remark}

\begin{document}

\title[]{Compact actions and uniqueness of the group measure space decomposition of II$_1$ factors}

\author[Adrian Ioana]{Adrian Ioana}
\address{Department of Mathematics, University of California, San Diego, La Jolla, CA
92093-0112, USA}
\email{aioana@ucsd.edu}

\begin{abstract}
We prove that any II$_1$ factor $L^{\infty}(X)\rtimes\Gamma$ arising from a compact, free, ergodic, probability measure preserving action of a countable group $\Gamma$ with positive first $\ell^2$-Betti number, has a unique group measure space Cartan subalgebra, up to unitary conjugacy.
\end{abstract}

\maketitle

 \section {Introduction}

The {\it group measure space construction} of Murray and von Neumann associates to every probability measure preserving (p.m.p.) action $\Gamma\curvearrowright (X,\mu)$ of a countable group $\Gamma$, a finite von Neumann algebra $M=L^{\infty}(X)\rtimes\Gamma$  \cite{MvN36}. If the action is (essentially)  free and ergodic, then $M$ is a II$_1$ factor and $A=L^{\infty}(X)$ is a {\it Cartan subalgebra}, i.e. a maximal abelian von Neumann subalgebra whose normalizer, $\mathcal N_{M}(A)=\{u\in\mathcal U(M)|uAu^*=A\}$, generates $M$, as a von Neumann algebra.

During the last decade, S. Popa's {\it deformation/rigidity} theory has led to spectacular progress in the classification of  group measure space II$_1$ factors (see  the surveys \cite{Po07},\cite{Va10a}).  In particular, several large families of group measure space II$_1$ factors $L^{\infty}(X)\rtimes\Gamma$  have been shown to have a unique Cartan subalgebra \cite{OP07},\cite{OP08},\cite{CS11},\cite{CSU11},\cite{PV11a},\cite{PV11b} or a unique group measure space Cartan subalgebra \cite{Pe09},\cite{PV09},\cite{Io10},\cite{FV10},\cite{IPV10},\cite{CP10}, \cite{HPV10},\cite{Va10b},\cite{Io11}. 
Such ``unique Cartan subalgebra" results are extremely useful because they enable one to reduce the classification of the factors $L^{\infty}(X)\rtimes\Gamma$,  up to isomorphism, to the classification of the corresponding actions,   up to {\it orbit equivalence}. 
Indeed, by \cite{Si55}, \cite{FM77},
 two free, ergodic actions $\Gamma\curvearrowright X$ and $\Lambda\curvearrowright Y$ are orbit equivalent if and only if
 there exists an isomorphism of group measure space factors $L^{\infty}(X)\rtimes\Gamma\cong L^{\infty}(Y)\rtimes\Lambda$ which  identifies the Cartan subalgebras $L^{\infty}(X)$ and $L^{\infty}(Y)$. Recall that the  actions $\Gamma\curvearrowright X$ and $\Lambda\curvearrowright Y$ are orbit equivalent if there exists an isomorphism of probability spaces $\theta:X\rightarrow Y$ such that $\theta(\Gamma\cdot x)=\Lambda\cdot\theta(x)$, almost everywhere;  for recent progress on orbit equivalence, see the surveys \cite{Fu09},\cite{Ga10}.

In their breakthrough work \cite{OP07}, S. Popa and N. Ozawa proved that II$_1$ factors $L^{\infty}(X)\rtimes\Gamma$ associated with free, ergodic, compact actions of  free groups $\Gamma=\mathbb F_n$ and their direct products $\Gamma=\Bbb F_{n_1}\times...\times \Bbb F_{n_k}$  have a unique Cartan subalgebra, up to unitary conjugacy. Their proof makes crucial use of the fact that free groups satisfy  a strong form of {\it weak amenability}, called the {\it complete metric approximation property} ({\it c.m.a.p.}). Subsequently, the results of \cite{OP07} have been extended to other classes of groups $\Gamma$: lattices in products of SO$(n,1)$  and SU$(n,1)$ \cite{OP08},  groups having positive first $\ell^2$-Betti number and the c.m.a.p. \cite{Da10},\cite{Si10},  hyperbolic groups \cite{CS11} and direct products of hyperbolic groups \cite{CSU11} (the last two results notably using  \cite{Oz07} and \cite{Oz10}). Note, however, that in all of these cases one requires that both the action $\Gamma\curvearrowright X$ is compact and, implicitely, that the group $\Gamma$ is weakly amenable.

In a very recent breakthrough, S. Popa and S. Vaes obtained wide-ranging unique Cartan subalgebra results by removing the compactness assumption on the action. More precisely, they showed that if $\Gamma$ is either a weakly amenable group with  $\beta_1^{(2)}(\Gamma)>0$  \cite{PV11a} or a hyperbolic group \cite{PV11b} (or a direct product of groups in one of these two classes), then II$_1$ factors $L^{\infty}(X)\rtimes\Gamma$ associated with {\it arbitrary} free, ergodic actions of $\Gamma$ have a unique Cartan subalgebra, up to unitary conjugacy. 

The main result of \cite{OP07} shows in particular that any II$_1$ factor arising from a free, ergodic, compact action of a free group $\Bbb F_n$ has a unique  group measure space Cartan subalgebra.  The original motivation for our paper was to give a proof of this result which does not rely on the weak amenability of $\Bbb F_n$. By combining our recent work \cite{Io11} with techniques from \cite{OP07} and \cite{CP10}, we were able to not only prove this, but to show the following:

\begin{theorem}\label{main}  Let  $\Gamma\curvearrowright (X,\mu)$ be a free, ergodic, compact, p.m.p. action of a countable group $\Gamma$ with  $\beta_1^{(2)}(\Gamma)>0$. 

Then the II$_1$ factor $M=L^{\infty}(X)\rtimes\Gamma$ has a unique group measure space Cartan subalgebra, up to unitary conjugacy. More precisely, if $M=L^{\infty}(Y)\rtimes\Lambda$, for any free, ergodic, p.m.p. action $\Lambda\curvearrowright (Y,\nu)$, then we can find a unitary $u\in M$ such that $uL^{\infty}(X)u^*=L^{\infty}(Y)$.
\end{theorem}

Recall that a p.m.p. action $\Gamma\curvearrowright (X,\mu)$ is called {\it compact} if we can decompose $L^2(X)=\bigoplus_{n}\mathcal H_n$, where $\{\mathcal H_n\}_ n$ is a sequence of $\Gamma$-invariant, finite dimensional subspaces of $L^2(X)$. For instance, if $G$ is a compact group containing a dense copy of $\Gamma$ and $m_G$ denotes the Haar measure of $G$, then the left multiplication action $\Gamma\curvearrowright (G,m_G)$ is a free, ergodic, p.m.p. compact action.

The condition $\beta_1^{(2)}(\Gamma)>0$ is equivalent to $\Gamma$ being non--amenable and having an unbounded cocycle into its regular representation $\lambda:\Gamma\rightarrow\mathcal U(\ell^2\Gamma)$ \cite{PT07}, and is satisfied by any free product group $\Gamma=\Gamma_1*\Gamma_2$ with $|\Gamma_1|\geqslant 2$ and $|\Gamma_2|\geqslant 3$. 

Note that in the case $\Gamma$ does not have Haagerup's property, Theorem \ref{main} is a consequence of the main result of \cite{Pe09}.

A far-reaching conjecture predicts that II$_1$ factors $M=L^{\infty}(X)\rtimes\Gamma$ arising from free, ergodic, p.m.p. actions of groups $\Gamma$ with $\beta_1^{(2)}(\Gamma)>0$ have a unique Cartan subalgebra. Recently, there has been some progress on the ``group measure space" version of the conjecture. Firstly, it  has been shown in  \cite{PV09} that if $\Gamma=\Gamma_1*\Gamma_2$, where $\Gamma_1$ is an infinite property (T) group and $\Gamma_2$ is a non-trivial group, then
$M$ has a unique group measure space Cartan subalgebra. More generally, the same  has been proven in \cite{CP10} (see also \cite{Va10b}) under the assumption that $\Gamma$ admits
a non-amenable subgroup with the relative property (T).  Lastly, we proved in \cite{Io11} that $M$ also has a unique group measure space Cartan subalgebra if the action $\Gamma\curvearrowright (X,\mu)$ is rigid.

As mentioned above,  this conjecture has been very recently established in full generality for arbitrary actions of weakly amenable groups $\Gamma$ with $\beta_1^{(2)}(\Gamma)>0$ \cite{PV11a}.

Theorem \ref{main} provides further positive evidence towards this conjecture.  It implies that any residually finite group $\Gamma$ satisfying $\beta_1^{(2)}(\Gamma)>0$ admits at least one action whose II$_1$ factor has a unique group measure space decomposition. Indeed, if $\{\Gamma_n\}_{n}$ is a descending chain of normal, finite index subgroups of $\Gamma$ with trivial intersection, then the left multiplication action $\Gamma\curvearrowright (G,m_G)$ on the {\it profinite completion}  $G=\varprojlim\Gamma/\Gamma_n$  is a free, ergodic, compact p.m.p. action.

To outline the proof of Theorem \ref{main}, denote $A=L^{\infty}(X)$ and $M=A\rtimes\Gamma$. Consider another group measure space decomposition $M=B\rtimes\Lambda$ and let $\Delta:M\rightarrow M\overline{\otimes}M$ be the associated {\it comultiplication} \cite{PV09}. 
In the first part of the proof, we show that a corner of $\Delta(A)$ embeds into $M\overline{\otimes}A$, in the sense of  \cite{Po03}.
To achieve this, we combine the fact that the action $\Delta(\Gamma)\curvearrowright \Delta(A)$ is compact (hence,  a fortiori, weakly compact) with a result from \cite{OP07} relating weak compactness of actions to uniform convergence of deformations, and results of \cite{CP10} on malleable deformations coming from cocycles into the regular representation (see Section 4).
In the second part of the proof, we use techniques from \cite{Io11} to conclude that if $\Gamma$ has positive first $\ell^2$-Betti number, then the condition $\Delta(A)\prec M\overline{\otimes}A$ automatically implies that $A$ and $B$ are conjugate (see Section 3).

\section {Preliminaries}

Throughout the paper we work with {\it tracial von Neumann algebras} $(M,\tau)$, i.e. von Neumann algebras $M$ endowed with a faithful normal tracial state $\tau$, and assume that all von Neumann algebras have separable predual.

\subsection {Intertwining-by-bimodules} We first recall from  \cite [Theorem 2.1 and Corollary 2.3]{Po03} S. Popa's powerful {\it intertwining-by-bimodules} technique.

\begin {theorem}\label{inter}\cite{Po03} Let $(M,\tau)$ be a tracial von Neumann algebra and $P,Q\subset M$ be von Neumann subalgebras. 
Then the following are equivalent:

\begin{itemize}

\item There exist  non-zero projections $p\in P, q\in Q$, a $*$-homomorphism $\phi:pPp\rightarrow qQq$  and a non-zero partial isometry $v\in qMp$ such that $\phi(x)v=vx$, for all $x\in pPp$.

\item There is no sequence $u_n\in\mathcal U(P)$ satisfying $||E_Q(xu_ny)||_2\rightarrow 0$, for all $x,y\in M$.
\end{itemize}

If one of these conditions holds true,  then we say that {\it a corner of $P$ embeds into $Q$ inside $M$} and write $P\prec_{M}Q$.
\end{theorem}

\subsection {Relative amenability}
Let $(M,\tau)$ be a tracial von Neumann algebra. Recall that $M$ is said to be {\it amenable} if we can find a sequence $\xi_n\in L^2(M)\overline{\otimes}L^2(M)$ such that $\langle x\xi_n,\xi_n\rangle\rightarrow\tau(x)$ and $||[x,\xi_n]||_2\rightarrow 0$, for every $x\in M.$ By A. Connes' celebrated theorem \cite {Co76} this is equivalent to $M$ being {\it approximately finite dimensional}.

Now, let $Q\subset M$ be a von Neumann subalgebra. {\it Jones' basic construction} $\langle M,e_Q\rangle$ is defined as the von Neumann subalgebra of $\mathbb B(L^2(M))$ generated by $M$ and the orthogonal projection $e_Q$ from $L^2(M)$ onto $L^2(Q)$.
Recall that $\langle M,e_Q\rangle$ is equipped with a faithful semi-finite  trace given by $Tr(xe_QyL)=\tau(xy)$ for all $x,y\in M$. We denote by $L^2(\langle M,e_Q\rangle)$ the associated Hilbert space and endow it with the natural $M$-bimodule structure. 

Following \cite[Definition 2.2]{OP07} we say that
 a (not necessarily unital) von Neumann subalgebra $P\subset M$ is {\it amenable relative to $Q$ inside $M$} if there exists a sequence $\xi_n\in L^2(\langle M,e_Q\rangle)$ such that $\langle x\xi_n,\xi_n\rangle\rightarrow \tau(x)$, for every $x\in M$, and $||[y,\xi_n]||_2\rightarrow 0$, for every $y\in P$.

\begin {lemma}\cite{IPV10}\label{ipv} Let $(B,\tau)$ be a tracial von Neumann algebra and $\Lambda\curvearrowright (B,\tau)$ be a trace preserving action. Denote $M=B\rtimes\Lambda$. Define the $*$-homomorphism $\Delta:M\rightarrow M\overline{\otimes}M$ by letting $\Delta(b)=b\otimes 1$, for all $b\in B$, and $\Delta(v_s)=v_s\otimes v_s$, for every $s\in\Lambda$.
Let $P,Q\subset M$ be von Neumann subalgebras such that $P$ has no amenable direct summand and $Q$ is amenable. 

Then there is no non-zero projection $p\in\Delta(P)'\cap M\overline{\otimes}M$ such that the von Neumann algebra $\Delta(P)p$ is amenable relative to $M\overline{\otimes}Q$.
\end{lemma}

{\it Proof.} Since $B$ and $Q$ are amenable, the $M$-$M\overline{\otimes}M$-bimodule $_{\Delta(M)}L^2(\langle M\overline{\otimes}M,e_{M\overline{\otimes}Q}\rangle)_{M\overline{\otimes}M}$ is weakly contained in the coarse $M$-$M\overline{\otimes}M$-bimodule, 
$L^2(M)\overline{\otimes}L^2(M\overline{\otimes}M)$.  Repeating the proof of  \cite[Proposition 7.2.{\it 4}]{IPV10} now gives the conclusion.
\hfill$\blacksquare$

\subsection {Weakly compact actions} We next recall N. Ozawa and S. Popa's notion of weakly compact actions.  
Let  $(P,\tau)$ be a tracial von Neumann algebra and denote by $\bar{P}=\{\bar{x}|x\in P\}$ the complex conjugate von Neumann algebra.

\vskip 0.05in
\begin{definition}\label{def}\cite[Definition 3.1]{OP07}. A trace preserving action $\Gamma\curvearrowright^{\sigma}(P,\tau)$ is {\it weakly compact} if we can find a sequence of unit vectors $\eta_n\in L^2(P\overline{\otimes}\bar{P})_{+}$ such that

\begin{itemize}
\item $||\eta_n-(v\otimes\bar{v})\eta_n||_2\rightarrow 0$, for every $v\in\mathcal U(P)$.

\item $||\eta_n-(\sigma_g\otimes\bar{\sigma}_g)(\eta_n)||_2\rightarrow 0$, for every $g\in\Gamma$.

\item $\langle (x\otimes 1)\eta_n,\eta_n\rangle=\tau(x)=\langle\eta_n,(1\otimes\bar{x})\eta_n\rangle$, for every $x\in P$ and every $n$. 
\end{itemize}
\end{definition}

Note that these conditions force $P$ to be amenable. 
On the other hand, by \cite[Proposition 3.2]{OP07},  if $P$ is amenable and $\Gamma\curvearrowright^{\sigma} P$ is a {\it compact} action (i.e. if $L^2(P)$ is the direct sum of finite dimensional $\sigma(\Gamma)$-invariant subspaces)  then $\sigma$ is weakly compact. 
\vskip 0.05in
N. Ozawa and S. Popa showed in \cite[Theorem 4.9]{OP07} that weakly compact actions can be used to deduce relative amenability. An obvious modification of their proof (in the case $p=1$)  gives the following:

\begin{theorem}\cite[Theorem 4.9]{OP07}\label{op}
 Let $(M,\tau)$ be a tracial von Neumann algebra and $Q_1,..,Q_k\subset M$ be von Neumann subalgebras. Assume that there are a tracial von Neumann algebra $\tilde M\supset M$ and automorphisms $\{\alpha_t\}_{t\in\Bbb R}$ of $\tilde M$ such that
\begin{itemize}
 \item $\lim_{t\rightarrow 0}||\alpha_t(x)-x||_2=0$, for all $x\in M$, and

\item $L^2(\tilde M)\ominus L^2(M)$ is isomorphic as an $M$-bimodule to a submodule of a multiple of $\bigoplus_{j=1}^k L^2(\langle M,e_{Q_j}\rangle)$.  
\end{itemize}

Let $P\subset M$ be a von Neumann subalgebra and define $c_t:=\inf_{v\in\mathcal U(P)}||E_M(\alpha_t(v))||_2$, for $t\in\Bbb R$.   Let $\mathcal G\subset\mathcal N_{M}(P)$ be a subgroup which acts weakly compactly on $P$ by conjugation and denote $N=\mathcal G''$. 

Then either $c_t\rightarrow 1$, as $t\rightarrow 0$,  or there exist $j\in\{1,..,k\}$ and a non-zero projection $p_j\in N'\cap M$ such that $Np_j$ is amenable relative to $Q_j$. 
\end{theorem}

\subsection {Deformations coming from group cocycles}\label{def}
Let  $\Gamma\curvearrowright (A,\tau)$ be a trace preserving action and set $M=A\rtimes\Gamma$. 
Let $\pi:\Gamma\rightarrow\mathcal O(H_{\Bbb R})$ be an orthogonal representation, where $H_{\Bbb R}$ is a separable real Hilbert space. Also, let  $b:\Gamma\rightarrow H_{\Bbb R}$ be a {\it cocycle}, i.e. a map satisfying the identity $b(gh)=b(g)+\pi(g)b(h)$, for all $g,h\in\Gamma$.

From this data, T. Sinclair constructed a {\it malleable deformation} in the sense of S. Popa, i.e. a tracial von Neumann algebra $\tilde M\supset M$ and a 1-parameter group of automorphisms $\{\alpha_t\}_{t\in\Bbb R}$ of $\tilde M$ such that $||\alpha_t(x)-x||_2\rightarrow 0$, for all $x\in\tilde M$  (see \cite[Section 3]{Si10} and \cite[Section 3.1]{Va10b}).

To recall this construction, let $(D,\tau_0)$ be the unique tracial von Neumann generated by unitary elements $\omega(\xi)$, $\xi\in H_{\Bbb R}$, subject to the relations 
 $\omega(\xi+\eta)=\omega(\xi)\omega(\eta)$, $\omega(\xi)^*=\omega(-\xi)$ and $\tau_0(\omega(\xi))=\exp(-||\xi||^2)$, for all $\xi,\eta\in H_{\Bbb R}$. 
 Consider the Gaussian action $\Gamma\curvearrowright^{\sigma} D$ which on the generating functions $\omega(\xi)$ is given by $\sigma_g(\omega(\xi))=\omega(\pi(g)(\xi))$. Finally, let 
$\Gamma\curvearrowright D\overline{\otimes}A$ be the diagonal action and define $\tilde M=(D\overline{\otimes}A)\rtimes\Gamma$.

Then the formula   $$\alpha_t(u_g)=(\omega(tb(g))\otimes 1)u_g,\;\;\text{for all}\;\;\; g\in\Gamma,\;\;\;\text{and}\;\;\; \alpha_t(x)=x,\;\;\text{for all}\;\; x\in D\overline{\otimes}A$$ gives the desired 1-parameter group of automorphisms $\{\alpha_t\}_{t\in\Bbb R}$ of $\tilde M$. 

Moreover, the formula $\beta(\omega(\xi))=\omega(-\xi)$, for all $\xi\in H_{\Bbb R}$, and $\beta(x)=x$, for all $x\in M$, defines an automorphism of $\tilde M$ which satisfies $\beta^2=id$ and $\beta\circ\alpha_t=\alpha_{-t}\circ\beta$, for all $t\in\Bbb R$.
As shown in \cite[Lemma 2.1]{Po06a} the existence of such a $\beta$ implies a ``transversality" property for $\alpha_t$:

\begin{lemma}\cite {Po06a}\label{trans} For every $x\in M$ and $t\in\Bbb R$ we have that $$||\alpha_{2t}(x)-x||_2\leqslant 2||\alpha_t(x)-E_M(\alpha_t(x))||_2.$$ \end{lemma}

Next we record a direct consequence of S. Popa's {\it spectral gap principle} \cite[Lemma 1.5]{Po06b}. 

\begin{lemma} \cite{Po06b}\label{spec} Assume that $\pi$ is weakly contained in the regular representation of $\Gamma$. Let $B\subset M$ be a von Neumann subalgebra with no amenable direct summand. 

Then $\alpha_t\rightarrow id$ in $||.||_2$ uniformly  on $(B'\cap M)_1$ (the unit ball of the commutant of $B$).
\end{lemma}

For a proof, see e.g. \cite[Lemma 2.2]{Io11}.

Now, let $B\subset M$ be a von Neumann subalgebra. 
J. Peterson \cite[Theorem 4.5]{Pe06} and I. Chifan and J. Peterson  \cite[Theorem 2.5]{CP10}  proved that if $\alpha_t\rightarrow id$ uniformly on $(B)_1$ and $B\nprec_{M}A$ then $\alpha_t\rightarrow id$ uniformly on  $\mathcal N_{M}(B)$.

\begin{theorem}\label{cp1}\cite{Pe06},\cite{CP10} Assume that $\pi$ is mixing. Let  $B\subset M$ be a von Neumann subalgebra. Suppose that $\alpha_t\rightarrow id$ in $||.||_2$ uniformly  on $(B)_1$ and that $B\nprec_{M}A$.
 Denote  $P=\mathcal N_M(B)''$. 

Then $\alpha_t\rightarrow id$ in $||.||_2$ uniformly on $(P)_1$. 
\end{theorem}

 Conversely, I. Chifan and J. Peterson proved in \cite[Theorem 3.2]{CP10} that if $B$ is abelian and $\alpha_t\rightarrow id$ uniformly on  a ``large" sequence  $u_n\in\mathcal N_{M}(B)$, then $\alpha_t\rightarrow id$ uniformly on $(B)_1$.

\begin{theorem}\label{cp2}\cite{CP10} Assume that $\pi$ is mixing. Let  $B\subset M$ be an abelian von Neumann subalgebra. 
Assume that  we can find a  sequence $u_n\in\mathcal N_{M}(B)$   such that  

\begin{itemize}
\item $\alpha_t\rightarrow id$ in $||.||_2$ uniformly  on $\{u_n\}_{n\geqslant 1}$  and

\item $||E_A(xu_ny)||_2\rightarrow 0$, for all $x,y\in M$.
\end{itemize}

Then $\alpha_t\rightarrow id$ in $||.||_2$ uniformly  on $(B)_1$.
\end{theorem}

Theorems \ref{cp1} and \ref{cp2} have been proven in \cite{Pe06} and \cite{CP10} using J. Peterson's technique of unbounded derivations \cite{Pe06}. For proofs using the  automorphisms $\{\alpha_t\}_{t\in\Bbb R}$, see S. Vaes's paper \cite[Theorems 3.9 and 4.1]{Va10b}.

\section {Conjugacy of group measure space Cartan subalgebras}

In this section we assemble together the main technical results of \cite{Io11} to prove a new conjugacy criterion for group measure space Cartan subalgebras of II$_1$ factors $L^{\infty}(X)\rtimes\Gamma$ arising from actions of groups with positive first $\ell^2$-Betti number.

\begin{theorem}\label{conjugacy} Let $\Gamma$ be a countable group satisfying $\beta_1^{(2)}(\Gamma)>0$. Let $\Gamma\curvearrowright (X,\mu)$ be a free ergodic p.m.p. action and set $M=L^{\infty}(X)\rtimes\Gamma$. 
Assume that $M=L^{\infty}(Y)\rtimes\Lambda$ for a free ergodic p.m.p. action $\Lambda\curvearrowright (Y,\nu)$.  
Let $\Delta:M\rightarrow M\overline{\otimes}M$ be the $*$-homomorphism given by $\Delta(b)=b\otimes 1$ for all $b\in L^{\infty}(Y)$ and $\Delta(v_g)=v_g\otimes v_g$ for all $g\in\Lambda$, where $\{v_g\}_{g\in\Lambda}\subset M$ denote the canonical unitaries. 

If $\Delta(L^{\infty}(X))\prec_{M\overline{\otimes}M}M\overline{\otimes}L^{\infty}(X)$, then $L^{\infty}(X)$ and $L^{\infty}(Y)$ are unitarily conjugate.
\end{theorem}

{\it Proof.} Denote $A=L^{\infty}(X)$ and $B=L^{\infty}(Y)$.  Since $\beta_1^{(2)}(\Gamma)>0$, by \cite[Theorem 4.2]{Io11}, in order to get the conclusion, it suffices to find an amenable subgroup $\Lambda_0<\Lambda$ such that $A\prec_{M}B\rtimes\Lambda_0$. Let us assume by contradiction that this is false.

Since $\Delta(A)\prec M\overline{\otimes}A$, the second condition in Theorem \ref{inter} holds true. Thus, we can find a finite set $F\subset (M\overline{\otimes}M)_1$  and $\delta>0$ such that
$\sum_{x\in F}||E_{M\overline{\otimes}A}(\Delta(a)x)||_2^2\geqslant 2\delta,$ for every  $a\in \mathcal U(A).$ By Kaplansky's theorem we may assume that $F\subset 1\otimes (M)_1$, i.e. we have

\begin{equation}\label{eq1}\sum_{x\in F}||E_{M\overline{\otimes}A}(\Delta(a)(1\otimes x))||_2^2\geqslant 2\delta,\hskip 0.05in\forall a\in\mathcal U(A)\end{equation}

For $g\in\Lambda$, we define $f(g)=\sum_{x\in F}||E_A(v_{g}x)||_2^2$. Then we have the following:

{\bf Claim.} Let $g_1,..,g_m,h_1,..,h_m\in\Lambda$ and $\Lambda_1,..,\Lambda_m<\Lambda$ be amenable subgroups. Put $S=\cup_{i=1}^mg_i\Lambda_ih_i$. Then there exists $g\in\Lambda\setminus S$ such that $f(g)\geqslant\delta$.

{\it Proof of the claim.}
Since $\Lambda_i$ is amenable by our assumption we have that $A\nprec_{M}B\rtimes\Lambda_i$, for all $i\in\{1,..,m\}$.
It follows from S. Popa's intertwining-by-bimodules technique (see Theorem \ref{inter} and e.g. \cite[Remark 1.2]{Io11}) that we can find $a\in\mathcal U(A)$ such that 

\begin{equation}\label{eq2}\sum_{i=1}^m||E_{B\rtimes\Lambda_i}(v_{g_i}^*av_{h_i}^*)||_2^2\leqslant \frac{\delta}{|F|}\end{equation}

If we write $a=\sum_{g\in\Lambda}b_gv_g$, where $b_g\in B$, then \ref{eq2} implies that $\sum_{g\in S}||b_g||_2^2\leqslant\frac{\delta}{|F|}.$
On the other hand, since $E_{M\overline{\otimes}A}(\Delta(a)(1\otimes x))=\sum_{g\in\Lambda}b_gv_g\otimes E_A(v_gx)$, for every $x\in M$, \ref{eq1} can be rewritten as $\sum_{g\in\Lambda}f(g)||b_g||_2^2\geqslant 2\delta$.
Since $f(g)\leqslant |F|$, for all $g\in\Lambda$, combining the last two inequalities yields that
 $\sum_{g\in\Lambda\setminus S}f(g)||b_g||_2^2\geqslant\delta$. Finally, since $\sum_{g\in\Lambda\setminus S}||b_g||_2^2\leqslant ||a||_2^2=1$, the claim follows.\hfill$\square$

The claim shows that the conclusion of \cite[Lemma 3.2]{Io11} holds true. The second part of the proof of \cite[Theorem 3.1]{Io11}, which only uses the conclusion of \cite[Lemma 3.2]{Io11},  provides a decreasing sequence $\{\Lambda_k\}_{k\geqslant 1}$ of non-amenable subgroups of $\Lambda$ such that $A\prec_{M}B\rtimes (\cup_{k\geqslant 1}C(\Lambda_k))$, where $C(\Lambda_k)$ denotes the centralizer of $\Lambda_k$ in $\Lambda$.

Since $\beta_1^{(2)}(\Gamma)>0$, there is an unbounded cocycle $b:\Gamma\rightarrow\ell^2_{\Bbb R}\Gamma$ for the regular representation \cite{PT07}.
Let $\tilde M\supset M$ and $\{\alpha_t\}_{t\in\Bbb R}\subset$ Aut$(\tilde M)$ be the malleable deformation constructed in Section \ref{def}.

Now, the group $\cup_{k\geqslant 1}C(\Lambda_k)$ must be non-amenable by our assumption. Thus, $C(\Lambda_k)$ is non-amenable for some $k\geqslant 1$ and therefore its von Neumann algebra, $L(C(\Lambda_k))$, has no amenable direct summand. Lemma \ref{spec} implies that $\alpha_t\rightarrow id$ in $||.||_2$   uniformly on $(L(\Lambda_k))_1$.

 Since $\Lambda_k$ is non-amenable and $A$ is abelian, we have that $L(\Lambda_k)\nprec_{M}A$. By \cite[Corollary 2.3]{Po03} can find a sequence $g_n\in\Lambda_k$ such that $||E_A(xv_{g_n}y)||_2\rightarrow 0$, for all $x,y\in M$. 
Because $v_{g_n}\in\mathcal N_{M}(B)$ and $\alpha_t\rightarrow id$ in $||.||_2$ uniformly on $\{v_{g_n}\}_{n\geqslant 1}$,  applying Theorem \ref{cp2} gives that $\alpha_t\rightarrow id$ in $||.||_2$ uniformly on 
$(B)_1$. 

Finally, if $B\prec_{M}A$, then by \cite[Theorem A.1]{Po01} we get that $A$ and $B$ are conjugate. Otherwise, 
 if $B\nprec_{M}A$, then Theorem \ref{cp1} implies that $\alpha_t\rightarrow id$ in $||.||_2$ uniformly on $(M)_1$. In particular, we can  find $t>0$ such that $\Re\tau(\alpha_t(u_g)u_g^*)\geqslant\exp(-\frac{1}{4})$, for all $g\in\Gamma$. Since $\tau(\alpha_t(u_g)u_g)=\exp(-t^2||b(g)||^2)$ we derive that $||b(g)||\leqslant\frac{1}{2t}$ for all $g\in\Gamma$. This shows that $b$ is bounded, a contradiction.\hfill$\blacksquare$

\section{Proof of Theorem \ref{main}}

Let $\Gamma\curvearrowright (X,\mu)$ be a free ergodic p.m.p. compact action. Denote $A=L^{\infty}(X)$ and $M=A\rtimes\Gamma$.   Assume that $\Gamma$ has positive first $\ell^2$--Betti number and let $b:\Gamma\rightarrow\ell^2_{\Bbb R}\Gamma$ be an unbounded cocycle into the regular representation of $\Gamma$ \cite{PT07}.

Let $\Lambda\curvearrowright (Y,\nu)$ be a free ergodic p.m.p. action such that $M=B\rtimes\Lambda$, where $B=L^{\infty}(Y)$. Let $\{v_s\}_{s\in\Lambda}\subset M$ be the canonical unitaries. Let $\Delta:M\rightarrow M\overline{\otimes}M$ be the $*$--homomorphism given by $\Delta(b)=b\otimes 1$ for $b\in B$ and $\Delta(v_s)=v_s\otimes v_s$ for $s\in\Lambda$.

Our goal is to show that $A$ and $B$ are unitarily conjugate. If we denote $\mathcal A=M\overline{\otimes}A$ and $\mathcal M=M\overline{\otimes}M$ then by Theorem \ref{conjugacy} it suffices to prove that $\Delta(A)\prec_{\mathcal M}\mathcal A$. Assume this to be false.

Towards a contradiction, let $\tilde M\supset M$ and $\{\alpha_t\}_{t\in\mathbb R}\subset$ Aut$(\tilde M)$ be the malleable deformation defined in Section \ref{def}. 
We let $\tilde{\mathcal M}=M\overline{\otimes}\tilde M$ and denote by  $\tilde\alpha_t$ the automorphism $id\otimes\alpha_t$ of $\tilde{\mathcal M}$ for every $t\in\mathbb R$. Note that $\tilde{\mathcal M}\supset \mathcal M$ and $\{\tilde\alpha_t\}_{t\in\mathbb R}\subset$ Aut$(\tilde{\mathcal M}$) is precisely the malleable deformation associated to the cocycle $b$ and the obvious crossed product decomposition $\mathcal M=\mathcal A\rtimes\Gamma$.

Since $b$ is a cocycle for the regular representation, it is easy to see that the $\mathcal M$-bimodule $L^2(\tilde{\mathcal M})\ominus L^2(\mathcal M)$ is isomorphic to  $\bigoplus_{i=1}^{\infty}L^2(\langle \mathcal M,e_{\mathcal A}\rangle)$ (see  \cite[Lemma 3.5]{Va10b}).
Since the action $\Gamma\curvearrowright A$ is compact, the action $\Delta(\Gamma)\curvearrowright \Delta(A)$ is compact and thus weakly compact. By applying Theorem \ref{op} to $P=\Delta(A)$ and $N=\Delta(M)$ we deduce that either (1)  $\Delta(M)p$ is amenable relative to $\mathcal A$ for some non-zero projection $p\in \Delta(M)'\cap \mathcal M$ or (2) $\inf_{v\in\mathcal U(\Delta(A))}||E_{\mathcal M}(\tilde\alpha_t(v))||_2\rightarrow 1$ as $t\rightarrow 0$. 

Since $M$ has no amenable direct summand, the first case is ruled out by Lemma \ref{ipv}.
 
In the second case, we get that $\sup_{v\in\mathcal U(\Delta(A))}||\tilde\alpha_t(v)-E_{\mathcal M}(\tilde\alpha_t(v))||_2\rightarrow 0$ as $t\rightarrow 0$. By Lemma \ref{trans} we deduce that $\tilde\alpha_t\rightarrow id$ in $||.||_2$ uniformly on $(\Delta(A))_1$. Since  $\Delta(A)\nprec_{\mathcal M}\mathcal A$ by assumption and the von Neumann algebra generated by the normalizer of $\Delta(A)$ in $\mathcal M$ contains $\Delta(M)$, Theorem \ref{cp1} gives that $\tilde\alpha_t\rightarrow id$  in $||.||_2$ uniformly on $(\Delta(M))_1$. 

Since $\tilde\alpha_t(\Delta(v_s))=v_s\otimes\alpha_t(v_s)$ for $s\in\Lambda$ we derive that $\alpha_t\rightarrow id$ in $||.||_2$ uniformly on $\{v_s\}_{s\in\Lambda}$.
Since $\Lambda$ is non--amenable and $A$ is abelian we can find a sequence $\{s_n\}_{n\geqslant 1}$ satisfying $||E_A(xv_{s_n}y)||_2\rightarrow 0$ for all $x,y\in M$. 
Continuing exactly as in the end of the proof of Theorem \ref{conjugacy} provides the desired contradiction.\hfill$\blacksquare$
\vskip 0.03in
\noindent
\begin{remark}
It is easy to adapt the proof of the main result to show that if $t>0$ and $M^t=L^{\infty}(Y)\rtimes\Lambda$, for some free ergodic p.m.p. action $\Lambda\curvearrowright (Y,\nu)$, then $L^{\infty}(Y)$ is unitarily conjugate to $L^{\infty}(X)^t$. It follows from D. Gaboriau's work on $\ell^2$-Betti numbers for equivalence relations \cite{Ga01} that $M$ has trivial fundamental group, $\mathcal F(M)=\{1\}$.
\end{remark}
\vskip 0.02in
\begin{remark}
Let $\mathbb Z_p$ be the ring of $p$-adic integers, for some prime $p$. Consider the profinite action $\Gamma:=\mathbb Z^2\rtimes$ SL$_2(\mathbb Z)\curvearrowright \varprojlim\mathbb Z_p^2$ and denote by $M$ the associated II$_1$ factor. By\cite [Corollary D]{OP08} and \cite[Section 5.5]{PV09}, $M$ has two non-conjugate group measure space Cartan subalgebras, $L^{\infty}(\varprojlim\mathbb Z_p^2)$ and $L(\mathbb Z^2)$. This shows that the assumption that $\beta_1^{(2)}(\Gamma)>0$ from the hypothesis of our main result cannot be be replaced by the assumption that $\Gamma$ is merely non-amenable.
\end{remark}


\begin{thebibliography}{}
\bibitem [Co76]{Co76} A. Connes: {\it Classification of injective factors}, Ann. of Math. (2) {\bf 104} (1976), 73--115.
\bibitem [CP10]{CP10} I. Chifan, J. Peterson: {\it Some unique group-measure space decomposition results}, 
 preprint arXiv:1010.5194.
\bibitem [CS11]{CS11} I. Chifan, T. Sinclair: {\it On the structural theory of II$_1$ factors of negatively curved groups}, preprint arXiv:1103.4299.
\bibitem [CSU11]{CSU11} I. Chifan, T. Sinclair, B. Udrea: {\it On the structural theory of II$_1$ factors of negatively curved groups, II. Actions by product groups}, preprint arXiv:1108.4200.
\bibitem [Da10]{Da10} Y. Dabrowski: {\it A non-commutative Path Space approach to stationary free Stochastic Differential Equations}, preprint arXiv:1006.4351.
\bibitem [FM77]{FM77} J. Feldman, C.C. Moore: {\it Ergodic equivalence relations, cohomology, and von Neumann algebras, II}, Trans. Amer. Math. Soc. {\bf 234} (1977), 325--359.
\bibitem [FV10]{FV10} P. Fima, S. Vaes: {\it HNN extensions and unique group measure space decomposition of II$_1$ factors}, preprint arXiv:1005.5002, to appear in Trans. Amer. Math. Soc.
\bibitem [Fu09]{Fu09} A. Furman: {\it A survey of Measured Group Theory}, Geometry, Rigidity, and Group Actions, 296--374, The University of Chicago Press, Chicago and London, 2011.
\bibitem [Ga01]{Ga01} D. Gaboriau: {\it Invariants $\ell^2$ de relations d'equivalence et de groupes}, Publ.
Math. Inst. Hautes \'Etudes Sci., {\bf 95} (2002), 93--150.
\bibitem [Ga10]{Ga10} D. Gaboriau: {\it Orbit Equivalence and Measured Group Theory},   Proceedings of the ICM (Hyderabad, India, 2010), Vol. III, Hindustan Book Agency (2010), 1501--1527. 
\bibitem [HPV10]{HPV10} C. Houdayer, S. Popa, S. Vaes: {\it A class of groups for which every action is W$^*$--superrigid}, preprint arXiv:1010.5077, to appear in Groups Geom. Dyn.
\bibitem [Io10]{Io10} A. Ioana: {\it W$^*$-superrigidity for Bernoulli actions of property (T) groups},  J. Amer. Math. Soc. {\bf 24} (2011), 1175--1226.
\bibitem [Io11]{Io11} A. Ioana: {\it Uniqueness of the group measure space decomposition for Popa's $\mathcal H\mathcal T$ factors}, preprint arXiv:1104.2913.
\bibitem [IPV10]{IPV10} A. Ioana, S. Popa, S. Vaes: {\it A class of superrigid group von Neumann algebras}, preprint arXiv:1007.1412.
\bibitem [MvN36]{MvN36} F. Murray, J. von Neumann: {\it On rings of operators,} Ann. of Math. {\bf 37} (1936), 116--229.
\bibitem [Oz07]{Oz07} N. Ozawa: {\it Weak amenability of hyperbolic groups,} Groups Geom. Dyn., {\bf 2} (2008), 271--280.
\bibitem [Oz10]{Oz10} N. Ozawa: {\it  Examples of groups which are not weakly amenable}, preprint arXiv:1012.0613, to appear in Kyoto J. Math.
\bibitem [OP07]{OP07} N. Ozawa, S. Popa: {\it On a class of II$_1$ factors with at most one Cartan subalgebra}, Ann. of Math. (2), {\bf 172} (2010), 713--749.
\bibitem [OP08]{OP08} N. Ozawa, S. Popa: {\it On a class of II$_1$ factors with at most one Cartan subalgebra, II}, Amer. J. Math., {\bf 132} (2010), 841--866. 
\bibitem [Pe06]{Pe06} J. Peterson: {\it $L^2$--rigidity in von Neumann algebras},  Invent. Math.  {\bf 175}  (2009),  no. 2, 417--433.
\bibitem [Pe09]{Pe09} J. Peterson: {\it Examples of group actions which are virtually W*-superrigid}, 
 preprint arXiv:1002.1745.
\bibitem [Po01]{Po01} S. Popa: {\it On a class of type II$_1$ factors with Betti numbers invariants}, Ann. of Math. {\bf 163} (2006), 809--889. 
\bibitem [Po03]{Po03} S. Popa: {\it Strong Rigidity of II$_1$ Factors Arising from Malleable Actions of w-Rigid Groups. I.}, Invent. Math. {\bf 165} (2006), 369--408.
\bibitem [Po06a]{Po06a} S. Popa: {\it On the superrigidity of malleable actions with spectral gap}, 
J. Amer. Math. Soc. {\bf 21} (2008), 981--1000. 
\bibitem [Po06b]{Po06b} S. Popa: {\it On Ozawa's Property for Free Group Factors}, Int. Math. Res. Notices (2007) Vol. 2007,
article ID rnm036.
\bibitem [Po07]{Po07} S. Popa: {\it Deformation and rigidity for group 	
actions and von Neumann algebras},  International Congress of Mathematicians. 
Vol. I,  445--477, 
Eur. Math. Soc., Z$\ddot{\text{u}}$rich, 2007.
\bibitem [PV09]{PV09} S. Popa, S. Vaes: {\it Group measure space decomposition of II$_1$ factors and} 
{\it W$^*$--superrigidity}, Invent. Math. {\bf 182} (2010), no. 2, 371--417.
\bibitem [PV11a]{PV11a} S. Popa, S. Vaes: {\it Unique Cartan decomposition for II$_1$ factors arising from arbitrary actions of free groups}, preprint arXiv:1111.6951.
\bibitem [PV11b]{PV11b} S. Popa, S. Vaes: {\it Unique Cartan decomposition for II$_1$ factors arising from arbitrary actions of hyperbolic groups}, in preparation.
\bibitem [PT07]{PT07} J. Peterson, A Thom: {\it Group cocycles and the ring of affiliated operators}, Invent. Math., {\bf 185} (2011) 561--592.
\bibitem [Si55]{Si55} I.M. Singer: {\it Automorphisms of finite factors}, Amer. J. Math. {\bf 77} (1955), 117--133.
 \bibitem [Si10]{Si10} T. Sinclair: {\it Strong solidity of group factors from lattices in SO(n,1) and SU(n,1)},  J. Funct. Anal. {\bf 260} (2011), 3209-3221.
\bibitem [Va10a]{Va10a} S. Vaes: {\it Rigidity for von Neumann algebras and their invariants}, Proceedings of the ICM (Hyderabad, India, 2010), Vol. III, Hindustan Book Agency (2010),  1624--1650.
\bibitem [Va10b]{Va10b} S Vaes: {\it One-cohomology and the uniqueness of the group
measure space decomposition of a II$_1$ factor}, preprint arXiv:1012.5377.
\end{thebibliography}
\enddocument